\documentclass[12pt,twoside]{amsart}
\usepackage{amsmath, amsthm, amscd, amsfonts, amssymb, graphicx}
\usepackage[bookmarksnumbered, plainpages]{hyperref}
\usepackage{marvosym}
\usepackage{geometry}
\makeatletter
\renewcommand\normalsize{%
\@setfontsize\normalsize\@xpt\@xiipt
\abovedisplayskip 5\p@ \@plus1\p@ \@minus6\p@
\abovedisplayshortskip \z@ \@plus3\p@
\belowdisplayshortskip 6\p@ \@plus3\p@ \@minus3\p@
\belowdisplayskip \abovedisplayskip
\let\@listi\@listI}
\makeatother

\newtheorem{thm}{Theorem}[section]

\newtheorem{lem}[thm]{Lemma}

\newtheorem{defn}[thm]{Definition}
\newtheorem{rem}[thm]{\bf{Remark}}

\newtheorem{pf}[thm]{\bf{Proof}}

\numberwithin{equation}{section}

\begin{document}

\title{\bf Left-invariant Ricci collineations associated to Yano connections on three-dimensional Lorentzian Lie groups}
\author{Yu Tao}
\thanks{{\scriptsize
\hskip -0.4 true cm \textit{2010 Mathematics Subject Classification:}
53C40; 53C42.
 }}

\maketitle

\begin{abstract}
In this paper, we classify Left-invariant Ricci collineations associated to Yano connections on three-dimensional Lorentzian Lie groups.

\end{abstract}

\vskip 0.2 true cm


\pagestyle{myheadings}
\markboth{\rightline\scriptsize }
{\leftline{\scriptsize Ricci collineation}}

\bigskip
\bigskip


\section{ Introduction}
\indent The concept of symmetry was first proposed by the ancient Greeks. Since the emergence of natural philosophy, the introduction of symmetry has provided support for the laws of physics and the essence of the cosmos. Two outstanding theoretical achievements of the 20th century, relativity and quantum mechanics, fundamentally involve the concept of symmetry. From the perspective of mathematics and physics, it is of great value to study symmetry in depth. After that, mathematicians simplified the Einstein equation and gave the classification of spacetimes according to the structure of the corresponding Lie algebra,which contains Ricci and curvature collineations ,among others(see\cite{A6},\cite{A7},\cite{A8},\cite{A9},\cite{A10},\cite{A11}).\\
\indent In \cite{A8},Michael and Pantelis had completed two tasks, the first of which was to put forward an path, which simplified the complexity of the calculation and lead to the Matter collineation and the Ricci collineation associated to a given Riemannian metric to the computation of Killing vectors, the second of which was to apply this path to determine the space-time metric of local rotational symmetry orthogonal to all hypersurfaces, which admit proper Ricci collineations and proper Matter collineations. In\cite{C1}, E.Calvino-Louzao, J.Seoane-Bascoy, M.E.Vazquez-Abal and R.Vazquez-Lorenzohe had proved that the concept of left-invariant Ricci collineations proposed to denote that  Ricci tensors  have a value of zero with a left-invariant Ricci collineation ~~$\xi$, having the mathematical meaning of a differential homomorphism of preserving Ricci tensors. To determine all left-invariant Ricci collineations on three-dimensional Lie groups \cite{C1}, it must inherit its symmetries where any homothetic vector field is a Ricci collineation.\\
\indent Yano connection can be used to study some new properties of manifolds and Lie groups, which makes the study of Yano connection have important geometric significance. Based on Yano 's work on Levi Civita connection, Etayo.F and Santanaria.R introduced Yano connection on manifolds with product structure or complex structure in \cite{A1}.In \cite{C3}, Wang proposed a product structure~~$(J\overline{e}_1=\overline{e}_1,J\overline{e}_2=\overline{e}_2~~and~~J\overline{e}_3=-\overline{e}_3)$~~and computed canonical connections and Kobayashi-Nomizu connections and their curvature on three-dimensional Lorentzian Lie groups with this product structure.In \cite{C2}, S.Azami applied this product structure to Yano connection as well, calculated its Lie derivative of curvature and metric on the three-dimensional Lorentzian Lie group with this product structure and introduced affine generalized Ricci solitons related to the Yano connection, and classified the left-invariant affine generalized Ricci solitons related to the Yano connection on the three-dimensional Lorentzian Lie group.In \cite{C5}, J.Miao,J.Yang and J.Guan completed the classification of three-dimensional Lorentzian Lie groups where Ricci tensors associated with Yano connections are Codazzi tensors.\\
\indent In this paper we continue to pay attention to Ricci collineation which is, ${\rm L}_{\xi} \overline{\rm Ric^{*}}=0$. our motivation is to classify all left-invariant Ricci collineations associated to Yano connections on three-dimensional Lorentzian Lie group. More precisely,the concept of Ricci collineation was extended to Yano connection on three-deimensional Lorentzian Lie group equiped with a product structure form \cite{C3}. We classify all left-invariant Ricci collineations on seven different three-dimensional Lie groups.\\
 \indent In section 2 of preliminary, we recall the definition of Yano connection and its corresponding symmetric tensor and Ricci collineation on the three-dimensional Lorentzian Lie group with product structure. At the same time, we directly reduce the amount of simple calculation and make some explanations in advance.\\
 \indent In section 3, we classify all left-invariant Ricci collineations associated to Yano connection on three-dimensional unimodular Lorentzian Lie groups. In section 4, we classify all left-invariant Ricci collineations associated to Yano connection on three-dimensional non-unimodular Lorentzian Lie groups.

\vskip 1 true cm
\section{Preliminaries}
\vskip 0.5 true cm
\noindent{\bf 2.1 The symmetric tensors associated to Yano connection }\\
\indent In this section, we recall the definition of Yano connection on a three-dimensional Lorentzian Lie group with product structure J and its corresponding symmetric tensor.\\
\begin{rem} $\nabla^{LC}$~~ is the Levi-Civita connection of $\{G_i\}_{i=1\cdot\cdot\cdot 7}$ and ~~$R^{LC}$ ~~is its curvature tensor, accompanied with the convention
\begin{equation}
R^{LC}(X,Y)Z=\nabla^{LC}_X\nabla^{LC}_YZ-\nabla^{LC}_Y\nabla^{LC}_XZ-\nabla^{LC}_{[X,Y]}Z.\notag
\end{equation}
\end{rem}
\begin{rem} The Ricci tensor of $(G_i,g)_{i=1\cdot \cdot \cdot 7} $ is defined by
\begin{equation}{\rm Ric}^{LC}(X,Y)=-g(R^{LC}(X,\overline{e}_1)Y,\overline{e}_1)-g(R^{LC}(X,\overline{e}_2)Y,\overline{e}_2)+g(R^{LC}(X,\overline{e}_3)Y,\overline{e}_3),\notag
\end{equation}

where $\{\overline{e}_1,\overline{e}_2,\overline{e}_3\}$ is a pseudo-orthonormal basis, with $\overline{e}_3$ timelike and the Ricci operator ~~$Ric^{LC}$~~ is given by
\begin{equation}{\rm Ric}^{LC}(X,Y)=g({\rm Ric}^{LC}(X),Y).\notag
\end{equation}
\end{rem}
\begin{rem}
A product structure J on $\{G_i\}_{i=1\cdot\cdot\cdot 7}$ was defined by
\begin{equation}J\overline{e}_1=\overline{e}_1,~J\overline{e}_2=\overline{e}_2,~J\overline{e}_3=-\overline{e}_3,\notag
\end{equation}
\end{rem}
then $J^2={\rm id}$ and $g(J\overline{e}_j,J\overline{e}_j)=g(\overline{e}_j,\overline{e}_j)$. \\
\begin{rem}
Let ~~$\nabla^{*}$~~ be the Yano connection by \cite{A5} following:
\begin{equation}\nabla^{*}_XY=\nabla^{LC}_XY-\frac{1}{2}(\nabla^{LC}_YJ)JX-\frac{1}{4}[(\nabla^{LC}_XJ)JY-(\nabla^{LC}_{JX}J)Y],\notag
\end{equation}
\end{rem}
\begin{rem}
Let ~~$R^{*}$~~ be its curvature tensors given by
\begin{align}\notag
R^{*}(X,Y)Z&=[\nabla^{*}_X,\nabla^{*}_Y]Z-\nabla^{*}_{[X,Y]}Z\\\notag
&=\nabla^{*}_X\nabla^{*}_YZ-\nabla^{*}_Y\nabla^{*}_XZ-\nabla^{*}_{[X,Y]}Z,\notag
\end{align}
\end{rem}
Then the Ricci tensor of ~~$(G_i,g,J)_{i=1\cdot\cdot\cdot 7}$ ~~associated to the Yano connection
is defined by
\begin{equation}{\rm Ric}^{*}(X,Y)=-g(R^{*}(X,\overline{e}_1)Y,\overline{e}_1)-g(R^{*}(X,\overline{e}_2)Y,\overline{e}_2)+g(R^{*}(X,\overline{e}_3)Y,\overline{e}_3),\notag
\end{equation}
 The Ricci operators ${\rm Ric}^{*}$ is given by
\begin{equation}{\rm Ric}^{*}(X,Y)=g({\rm Ric}^{*}(X),Y).\notag
\end{equation}
\begin{rem}
Let ~~$\overline{\rm Ric^{*}}$~~ be symmetric tensors associated Yano connection
\begin{equation}\overline{\rm Ric^{*}}(X,Y)=\frac{{\rm Ric}^{*}(X,Y)+{\rm Ric}^{*}(Y,X)}{2}\notag
\end{equation}
and
\begin{equation}\overline{\rm Ric^{*}}(X,Y)=g(\overline{\rm Ric^{*}}(X),Y).\notag
\end{equation}
\end{rem}
According above analysis so far we can prove  that $\overline{\rm Ric^{*}}~~$ are symmetric tensors.\\
There are different and specific three-dimensional Lorentzian Lie groups~~$(G_i,g,J)_{i=1\cdot\cdot\cdot 7}$~~in \cite{A3,A4}(Theorem 2.1 and Theorem 2.2 \cite{A2}), and then we have complete specific classification about three-dimensional Lorentzian Lie groups . We shall denote the connected three-dimensional unimodular  and non-unimodular Lie groups, and having a corresponding Lie algebra $\{\mathfrak {g}\}_{i=1,\cdots,7}$~~by \cite{A4}.
\vskip 0.5 true cm
\noindent{\bf 2.2 Three-dimensional unimodular Lorentzian Lie groups }\\
\begin{thm}
Assume~ $(G,g,J)$~was a three-dimensional unimodular Lorentzian Lie group,equipped with a left-invariant Lorentzian metric~$\boldsymbol g$~and a product structure~$\boldsymbol J$~,and had a Lie algebra ~$\{\mathfrak {g}\}$~ being one of the following:\\
$\mathfrak {g}_1:$
\begin{equation}
  [\overline{e}_1,\overline{e}_2]=m \overline{e}_1-n \overline{e}_3,~~[\overline{e}_1,\overline{e}_3]=-me_1-n \overline{e}_2,~~[\overline{e}_2,\overline{e}_3]=n \overline{e}_1+m \overline{e}_2+m \overline{e}_3.\notag
\end{equation}
$\mathfrak {g}_2:$
\begin{equation}
  [\overline{e}_1,\overline{e}_2]=n \overline{e}_2-u \overline{e}_3,~~[\overline{e}_1,\overline{e}_3]=-u \overline{e}_2-ma \overline{e}_3,~~[\overline{e}_2,\overline{e}_3]=m \overline{e}_2,~~n\neq0.\notag
\end{equation}
$\mathfrak {g}_3:$
\begin{equation}
  [\overline{e}_1,\overline{e}_2]=-ue_3,~~[\overline{e}_1,\overline{e}_3]= -n \overline{e}_2,~~[\overline{e}_2,\overline{e}_3]=m \overline{e}_1~~.\notag
\end{equation}
$\mathfrak {g}_4:$
\begin{equation}
  [\overline{e}_1,\overline{e}_2]=-\overline{e}_2+(2v-n)\overline{e}_3,v=\pm1,~~[\overline{e}_1,\overline{e}_3]=-n \overline{e}_2+\overline{e}_3,~~[\overline{e}_2,\overline{e}_3]=m \overline{e}_1,~~.\notag
\end{equation}
\end{thm}
\vskip 0.1 true cm
\noindent{\bf 2.3 Three-dimensional non-unimodular Lorentzian Lie groups }\\
\begin{thm}
Assume~ $(G,g,J)$~was a three-dimensional non-unimodular Lorentzian Lie group,equipped with a left-invariant Lorentzian metric~$\boldsymbol g$~and a product structure~$\boldsymbol J$~,and had a Lie algebra ~$\{\mathfrak {g}\}$~ being one of the following:\\
$\mathfrak {g}_5:$
\begin{equation}
  [\overline{e}_1,\overline{e}_2]=0,~~[\overline{e}_1,\overline{e}_3]=m \overline{e}_1+n \overline{e}_2,~~[\overline{e}_2,\overline{e}_3]=u \overline{e}_1+v \overline{e}_2,~~m+v\neq0~~,mu+nv=0.\notag
\end{equation}
$\mathfrak {g}_6:$
\begin{equation}
  [\overline{e}_1,\overline{e}_2]=m \overline{e}_2+n \overline{e}_3,~~[\overline{e}_1,\overline{e}_3]=u \overline{e}_2+v \overline{e}_3,~~[\overline{e}_2,\overline{e}_3]=0,~~m+v\neq0~~,mu-nv=0.\notag
\end{equation}
$\mathfrak {g}_7:$
\begin{equation}
  [\overline{e}_1,\overline{e}_2]=-m \overline{e}_1-n \overline{e}_2-n \overline{e}_3,~~[\overline{e}_1,\overline{e}_3]=m \overline{e}_1+n \overline{e}_2+n \overline{e}_3,~~[\overline{e}_2,\overline{e}_3]=u \overline{e}_1+v \overline{e}_2+v \overline{e}_3,~~m+v\neq0~~,mu=0.\notag
\end{equation}
\end{thm}
\begin{defn}\label{2.4}
\vskip 0.5 true cm
A Lorentzian Lie group of $\{G_i\}_{i=1,\cdots,7}$ admits left-Ricci collineations if and only if it satisfies:\\
\begin{equation}
  {\rm L}_{\xi}\overline{\rm Ric^{*}}=0,
\end{equation}
where ~~$\xi~~$ is a element of three-dimensional Lie group and~~ $\xi=\lambda_1 \overline{e}_1+\lambda_2 \overline{e}_2+\lambda_3 \overline{e}_3$~~for constant ~~$\lambda_1,\lambda_2,\lambda_3 \in\mathbb{R}$,~~and ~~$\{\overline{e}_1,\overline{e}_2,\overline{e}_3\}$ ~~is a pseudo-orthonormal basis with $\overline{e}_3$ timelike.
Let
\begin{equation}{\rm L}_{\xi} \overline{\rm Ric^{*}}(X,Y)=\xi(\overline{\rm Ric^{*}} (X,Y)-\overline{\rm Ric^{*}}([\xi,X],Y)-\overline{\rm Ric^{*}}(X,[\xi,Y]),
\end{equation}
\end{defn}
\indent According to this definition,we have
\begin{align}\notag
{\rm L}_{\overline{e}_1} \overline{\rm Ric^{*}}(\overline{e}_1,\overline{e}_1)=\overline{e}_1(\overline{\rm Ric^{*}}(\overline{e}_1,\overline{e}_1)-\overline{\rm Ric^{*}}([\overline{e}_1,\overline{e}_1],\overline{e}_1)-\overline{\rm Ric^{*}}(\overline{e}_1,[\overline{e}_1,\overline{e}_1])=0,~~\\\notag
{\rm L}_{\overline{e}_2} \overline{\rm Ric^{*}}(\overline{e}_2,\overline{e}_2)=\overline{e}_2\overline{\rm Ric^{*}} (\overline{e}_2,\overline{e}_2)-\overline{\rm Ric^{*}}([\overline{e}_2,\overline{e}_2],\overline{e}_2)-\overline{\rm Ric^{*}}(\overline{e}_2,[\overline{e}_2,\overline{e}_2])=0,~~\\\notag
{\rm L}_{\overline{e}_3} \overline{\rm Ric^{*}}(\overline{e}_3,\overline{e}_3)=\overline{e}_3(\overline{\rm Ric^{*}} (\overline{e}_3,\overline{e}_3)-\overline{\rm Ric^{*}}([\overline{e}_3,\overline{e}_3],\overline{e}_3)-\overline{\rm Ric^{*}}(\overline{e}_3,[\overline{e}_3,\overline{e}_3])=0,~~\\\notag
\end{align}
in the following discussion,we will not mention the calculation of these six items.\\
Owing to $\overline{\rm Ric^{*}}~~$ are symmetric tensors,we have
\begin{align}
{\rm L}_{\xi} \overline{\rm Ric^{*}} (\overline{e}_1,\overline{e}_2)={\rm L}_{\xi} \overline{\rm Ric^{*}} (\overline{e}_2,\overline{e}_1) ,~~\notag
{\rm L}_{\xi}\overline{\rm Ric^{*}}(\overline{e}_1,\overline{e}_3)={\rm L}_{\xi} \overline{\rm Ric^{*}} (\overline{e}_3,\overline{e}_1) ,~~\notag
\end{align}
\begin{rem}
 $\mathbb{V}_{\mathbb{R}\mathbb{C}}$ is spanned by Ricci collineations.
\end{rem}
\section{ Left-invariant Ricci collineations associated to Yano connections on three-dimensional unimodular Lorentzian Lie groups}
\vskip 0.5 true cm
\noindent{\bf 3.1 Left-invariant Ricci collineation of $G_1$}\\
\vskip 0.5 true cm
By (30) in \cite{C2},~~we have
\\
\begin{lem}
Ricci symmetric tensors of $(G_1,g,J)$  associated to Yano connection ~~$\nabla^{*}$~~ are given by
\begin{align}
\overline{\rm Ric^{*}}(\overline{e}_i,\overline{e}_j)=\left(\begin{array}{ccc}
-(m^2+n^2)&mn&-\frac{1}{2}mn\\
mn&-(m^2+n^2)&\frac{1}{2}m^2\\
-\frac{1}{2}mn&\frac{1}{2}m^2&0
\end{array}\right).
\end{align}
\end{lem}
By (2.2) and (3.1), we have
\begin{lem}
For $G_1$, the following equalities hold
\begin{align}
&{\rm L}_{\overline{e}_1}\overline{\rm Ric^{*}}(\overline{e}_1,\overline{e}_2)=m(m^2+\frac{1}{2} n^2),~~{\rm L}_{\overline{e}_1} \overline{\rm Ric^{*}}(\overline{e}_1,\overline{e}_3)=-m^3,~~{\rm L}_{\overline{e}_1}\overline{\rm Ric^{*}}(\overline{e}_2,\overline{e}_2)=-m^2n,~~\\\notag
&{\rm L}_{\overline{e}_1}\overline{\rm Ric^{*}}(\overline{e}_2,\overline{e}_3)=n(\frac{1}{2}m^2-n^2),~~{\rm L}_{\overline{e}_1} \overline{\rm Ric^{*}}(\overline{e}_3,\overline{e}_3)=0,~~{\rm L}_{\overline{e}_2} \overline{\rm Ric^{*}}(\overline{e}_1,\overline{e}_1)=-m(2m^2+n^2)~~\\\notag
&{\rm L}_{\overline{e}_2}\overline{\rm Ric^{*}}(\overline{e}_1,\overline{e}_2)=\frac{1}{2}m^2 n,~~{\rm L}_{\overline{e}_2}\overline{\rm Ric^{*}}(\overline{e}_1,\overline{e}_3)=n^3~~,{\rm L}_{\overline{e}_2} \overline{\rm Ric^{*}}(\overline{e}_2,\overline{e}_3)=\frac{1}{2}m^3,\\\notag
&{\rm L}_{\overline{e}_2}\overline{\rm Ric^{*}}(\overline{e}_3,\overline{e}_3)=m(n^2-m^2)~~,{\rm L}_{\overline{e}_3} \overline{\rm Ric^{*}}(\overline{e}_1,\overline{e}_1)=2m^3,~~{\rm L}_{\overline{e}_3} \overline{\rm Ric^{*}}(\overline{e}_1,\overline{e}_2)=-\frac{1}{2}m^2n~~\\\notag
&{\rm L}_{\overline{e}_3}\overline{\rm Ric^{*}}(\overline{e}_1,\overline{e}_3)=0,{\rm L}_{\overline{e}_3}\overline{\rm Ric^{*}}(\overline{e}_2,\overline{e}_2)=-m^3~~,{\rm L}_{\overline{e}_3} \overline{\rm Ric^{*}}(\overline{e}_2,\overline{e}_3)=\frac{1}{2}m(m^2-n^2).~~\\\notag
\end{align}
\end{lem}
By (3.2) and definition 2.9, we have
\vskip 0.5 true cm
\begin{thm}
the unimodular Lorentzian Lie group $(G_1,g,J)$ does not admit left-invariant Ricci collineations associated to Yano connection $\nabla^{*}$ .
\end{thm}
\begin{pf}
If the unimodular Lorentzian Lie group $(G_1,g,J)$ admits left-invariant Ricci collineations associated to Yano connection $\nabla^{*}$, then
${\rm L}_{\xi} \overline{\rm Ric_{ij}^{*}}=0,for ~~i=1,2,3 ~~and~~ j=1,2,3$,~~ so
\vskip 0.5 true cm
\begin{align}
\left\{\begin{array}{l}
-m(2m^2+n^2)\lambda_2+2m^3\lambda_3=0,\\
\\
m(m^2+\frac{1}{2}n^2)\lambda_1+\frac{1}{2}m^2 n\lambda_2-\frac{1}{2}m^2n\lambda_3=0,\\
\\
-m^3\lambda_1+n^3\lambda_2=0,\\
\\
m^2n\lambda_1+m^3\lambda_3=0,\\
\\
n(\frac{1}{2}m^2-n^2)\lambda_1+\frac{1}{2}m^3\lambda_2+\frac{1}{2}m(m^2-n^2)\lambda_3=0,\\
\\
m(n^2-m^2)\lambda_2=0.\\
\end{array}\right.
\end{align}
\vskip 0.5 true cm
  We analyze each one of these factors by separate. Because $m\neq 0,$~~then we have
  \vskip 0.5 true cm
 \begin{align}
\left\{\begin{array}{l}
-(2m^2+n^2)\lambda_2+2m^2\lambda_3=0,\\
\\
(n-m)(n+m)\lambda_2=0,\\
\\
n\lambda_1+m\lambda_3=0,\\
\\
-m^3\lambda_1+n^3\lambda_2=0,\\
\\
(m^2+\frac{1}{2}n^2)\lambda_1+\frac{1}{2}mn(\lambda_2-\lambda_3)=0,\\
\\
-n(\frac{1}{2}m^2+n^2)\lambda_1+\frac{1}{2}m^3\lambda_2+\frac{1}{2}m(m^2-n^2)\lambda_3=0.\\
\end{array}\right.
\end{align}
\vskip 0.5 true cm
Now if $m=n,$~~we have~
\begin{align}
\left\{\begin{array}{l}
m(\lambda_1+\lambda_3)=0,\\
\\
m^3(\lambda_2-\lambda_1)=0.\\
\end{array}\right.
\end{align}
As $~~m\neq0$,~~naturelly we have ~~$\lambda_1=\lambda_2=\lambda_3.$~~ Put this into the first equation of system (3.4) and get they are equal to zero, So~~$\mathbb{V}_{\mathbb{R}\mathbb{C}}$~~ is empty set.\\
And if $m+n=0,$~~we have~
\begin{align}
\left\{\begin{array}{l}
m(\lambda_1-\lambda_3)=0,\\
\\
m^3(\lambda_1+\lambda_2)=0.\\
\end{array}\right.
\end{align}
As $~~m\neq0$,~~also we have ~~$\lambda_1=\lambda_2=\lambda_3=0.$~~ So~~$\mathbb{V}_{\mathbb{R}\mathbb{C}}$~~ is empty set as well.\\
In summary, the unimodular Lorentzian Lie group $(G_1,g,J)$ does not admit left-invariant Ricci collineations associated to the Yano connection $\nabla^{*}$ .
\end{pf}
\vskip 0.5 true cm
\vskip 0.5 true cm
\noindent{\bf 3.2 Left-invariant Ricci collineation of $G_2$}\\
\vskip 0.5 true cm
By (19) in \cite{C2}, we have
\\
\begin{lem}
Ricci symmetric tensors of $(G_2,g,J)$  associated to Yano connection~~$\nabla^{*}$~~ are given by
\begin{align}
\overline{\rm Ric^{*}}(\overline{e}_i,\overline{e}_j)=\left(\begin{array}{ccc}
-(b^2+u^2)&0&0\\
0&-\left(n^2+mu\right)&-\frac{1}{2}mn\\
0&-\frac{1}{2}mn&0
\end{array}\right).
\end{align}
\end{lem}
By (2.2) and (3.7), we have
\begin{lem}
For $G_2$, the following equalities hold
\begin{align}
&{\rm L}_{\overline{e}_1} \overline{\rm Ric^{*}}(\overline{e}_1,\overline{e}_2)=0,~~{\rm L}_{\overline{e}_1}\overline{\rm Ric^{*}}(\overline{e}_1,\overline{e}_3)=0~~,{\rm L}_{\overline{e}_1} \overline{\rm Ric^{*}}(\overline{e}_2,\overline{e}_2)=n(2n^2+mu),\\\notag
&{\rm L}_{\overline{e}_1} \overline{\rm Ric^{*}}(\overline{e}_2,\overline{e}_3)=-u(n^2+mu)~~,{\rm L}_{\overline{e}_1} \overline{\rm Ric^{*}}(\overline{e}_3,\overline{e}_3)=-mnu,~~{\rm L}_{\overline{e}_2} \overline{\rm Ric^{*}}(\overline{e}_1,\overline{e}_1)=0~~\\\notag
&{\rm L}_{\overline{e}_2} \overline{\rm Ric^{*}}(\overline{e}_1,\overline{e}_2)=-n(n^2+\frac{1}{2}mu),~~{\rm L}_{\overline{e}_2} \overline{\rm Ric^{*}}(\overline{e}_1,\overline{e}_3)=m(u^2+\frac{1}{2}n^2)~~,{\rm L}_{\overline{e}_2} \overline{\rm Ric^{*}}(\overline{e}_2,\overline{e}_3)=0,\\\notag
&{\rm L}_{\overline{e}_2}\overline{\rm Ric^{*}}(\overline{e}_3,\overline{e}_3)=0~~,{\rm L}_{\overline{e}_3} \overline{\rm Ric^{*}}(\overline{e}_1,\overline{e}_1)=0,~~{\rm L}_{\overline{e}_3} \overline{\rm Ric^{*}}(\overline{e}_1,\overline{e}_2)=n^2(u-\frac{1}{2}m)~~\\\notag
&{\rm L}_{\overline{e}_3}\overline{\rm Ric^{*}}(\overline{e}_1,\overline{e}_3)=\frac{1}{2}mnu,~~{\rm L}_{\overline{e}_3}\overline{\rm Ric^{*}}(\overline{e}_2,\overline{e}_2)=0~~,{\rm L}_{\overline{e}_3} \overline{\rm Ric^{*}}(\overline{e}_2,\overline{e}_3)=0.~~\\\notag
\end{align}
\end{lem}
By (3.8) and definition2.9, we get
\vskip 0.5 true cm
\begin{thm}
the unimodular Lorentzian Lie group $(G_2,g,J)$ admits left-invariant Ricci collineations associated to Yano connection $\nabla^{*}$ if and only if\\
$(1)m=0,n\neq0,\lambda_1=0,n\lambda_2=u\lambda_3\in \mathbb{R},\mathbb{V}_{\mathbb{R}\mathbb{C}}=\langle \frac{u}{n}\overline{e}_2+\overline{e}_3\rangle$,\\
\\
$(2)m\neq0,n\neq0,u=\frac{1}{4}m,\lambda_1=0,\lambda_2=-\frac{nu}{n^2+2u^2}\lambda_3,\mathbb{V}_{\mathbb{R}\mathbb{C}}=\langle -\frac{nu}{n^2+2u^2}\overline{e}_2+\overline{e}_3\rangle$.\\
\end{thm}
\begin{pf}
If the unimodular Lorentzian Lie group $(G_2,g,J)$ admits left-invariant Ricci collineations associated to Yano connection $\nabla^{*}$, then
${\rm L}_{\xi} \widetilde{\rm Ric}_{ij}^0=0,for ~~i=1,2,3~~ and~~ j=1,2,3$,~~ so
\begin{align}
\left\{\begin{array}{l}
-n(n^2+\frac{1}{2}mu)\lambda_2+n^2(u-\frac{1}{2}m)\lambda_3=0,\\
\\
m(u^2+\frac{1}{2}n^2)\lambda_2+\frac{1}{2}mnu\lambda_3=0,\\
\\
n(2n^2+mu)\lambda_1=0,\\
\\
u(n^2+mu)\lambda_1=0,\\
\\
mnu\lambda_1=0.\\
\end{array}\right.
\end{align}
  We analyze each one of these factors by separate. \\
  Because $n\neq 0,$~~we assume that~~$m=0$~~firstly and the system (3.9) would reduce to
 \begin{align}
\left\{\begin{array}{l}
-n^3\lambda_2+n^2u\lambda_3=0,\\
\\
n^3\lambda_1=0,\\
\\
n^2u\lambda_1=0.\\
\end{array}\right.
\end{align}
Solving these equtions and having $\lambda_1=0,n\lambda_2=u\lambda_3$~~. Case~~ $(i)$ ~~is true.\\
Suppose~~$m\neq0$~~(whether u equals to zero or not, we stiil have ~~$\lambda_1=0$~~) and get
\begin{align}
\left\{\begin{array}{l}
(u^2+\frac{1}{2}n^2)\lambda_2+\frac{1}{2}nu\lambda_3=0,\\
\\
n(n^2+\frac{1}{2}mu)\lambda_2-n^2(u-\frac{1}{2}m)\lambda_3=0.\\
\end{array}\right.
\end{align} Then we assume that
\begin{align}
A=(a_{ij})
,~~\left\{\begin{array}{l}
a_{11}=u^2+\frac{1}{2}n^2,\\
\\
a_{12}=\frac{1}{2}nu,\\
\\
a_{21}=n(n^2+\frac{1}{2}mu),\\
\\
a_{22}=-n^2(u-\frac{1}{2}m).\\
\end{array}\right.~~
\end{align}
We get$$
\left|A\right|
=n^2(n^2+u^2)(u-\frac{1}{4}m)
$$
If $u=\frac{1}{4}m$,~~we have~~$
\left|A\right|
=0$.~~Put this into ~~$(3.11)$~~and get $\lambda_2=-\frac{nu}{n^2+2u^2}\lambda_3$,~~so~~$\mathbb{V}_{\mathbb{R}\mathbb{C}}=\langle -\frac{nu}{n^2+2u^2}\overline{e}_2+\overline{e}_3\rangle $~~ and case (i) is ture.\\
Notice that if $u\neq \frac{1}{4}m$,~~we have$
\left|A\right|
\neq0
$, from the relationship between Matrix and the solution of equation, then~~$\lambda_2=\lambda_3=0.$~~~~$\mathbb{V}_{\mathbb{R}\mathbb{C}}$~~is empty set and we complete the proof of theorem.\\
\end{pf}
\vskip 0.5 true cm
\noindent{\bf 3.3 Left-invariant Ricci collineation of $G_3$}\\
\vskip 0.5 true cm
By (8) in \cite{C2}, we have
\\
\begin{lem}
Ricci symmetric tensors of $(G_3,g,J)$  associated to Yano connection~~$\nabla^{*}$~~ are given by
\begin{align}
\overline{\rm Ric^{*}}(\overline{e}_i,\overline{e}_j)=\left(\begin{array}{ccc}
-nu&0&0\\
0&-mu&0\\
0&0&0
\end{array}\right).
\end{align}
\end{lem}
By (2.2) and (3.13), we have
\begin{lem}
For $G_3$, the following equalities hold
\begin{align}
&{\rm L}_{\overline{e}_1} \overline{\rm Ric^{*}}(\overline{e}_1,\overline{e}_2)=0,~~{\rm L}_{\overline{e}_1}\overline{\rm Ric^{*}}(\overline{e}_1,\overline{e}_3)=0~~,{\rm L}_{\overline{e}_1} \overline{\rm Ric^{*}}(\overline{e}_2,\overline{e}_2)=0,\\\notag
&{\rm L}_{\overline{e}_1} \overline{\rm Ric^{*}}(\overline{e}_2,\overline{e}_3)=-mnu~~,{\rm L}_{\overline{e}_1} \overline{\rm Ric^{*}}(\overline{e}_3,\overline{e}_3)=0,~~{\rm L}_{\overline{e}_2} \overline{\rm Ric^{*}}(\overline{e}_1,\overline{e}_1)=0~~\\\notag
&{\rm L}_{\overline{e}_2} \overline{\rm Ric^{*}}(\overline{e}_1,\overline{e}_2)=0,~~{\rm L}_{\overline{e}_2} \overline{\rm Ric^{*}}(\overline{e}_1,\overline{e}_3)=mnu~~,{\rm L}_{\overline{e}_2} \overline{\rm Ric^{*}}(\overline{e}_2,\overline{e}_3)=0,\\\notag
&{\rm L}_{\overline{e}_2}\overline{\rm Ric^{*}}(\overline{e}_3,\overline{e}_3)=0~~,{\rm L}_{\overline{e}_3} \overline{\rm Ric^{*}}(\overline{e}_1,\overline{e}_1)=0,~~{\rm L}_{\overline{e}_3} \overline{\rm Ric^{*}}(\overline{e}_1,\overline{e}_2)=0~~\\\notag
&{\rm L}_{\overline{e}_3}\overline{\rm Ric^{*}}(\overline{e}_1,\overline{e}_3)=0,~~{\rm L}_{\overline{e}_3}\overline{\rm Ric^{*}}(\overline{e}_2,\overline{e}_2)=0~~,{\rm L}_{\overline{e}_3} \overline{\rm Ric^{*}}(\overline{e}_2,\overline{e}_3)=0.~~\\\notag
\end{align}
\end{lem}
By (3.14) and definition 2.9, we have
\vskip 0.5 true cm
\begin{thm}
the unimodular Lorentzian Lie group $(G_3,g,J)$ admits left-invariant Ricci collineations associated to Yano connection $\nabla^{*}$ if and only if\\
\\
$(1)mnu=0,\mathbb{V}_{\mathbb{R}\mathbb{C}}=\langle \overline{e}_1,\overline{e}_2,\overline{e}_3\rangle,$\\
\\
$(2)mnu\neq0,\mathbb{V}_{\mathbb{R}\mathbb{C}}=\langle \overline{e}_3\rangle.$\\
\end{thm}
\begin{pf}
If the unimodular Lorentzian Lie group $(G_3,g,J)$ admits left-invariant Ricci collineations associated to the Yano connection $\nabla^{*}$, then
${\rm L}_{\xi} \overline{\rm Ric^{*}_{ij}}^0=0,i=1,2,3~~ and~~ j=1,2,3$, so
\begin{align}
\left\{\begin{array}{l}
mnu\lambda_1=0,\\
\\
mnu\lambda_2=0.\\
\end{array}\right.
\end{align}
  We analyze each one of these factors by separate. There is obvious that the value of~~ $\lambda_1$~~and ~~$\lambda_2$~~is related to the positive or negative of their coefficient. So when~~ $mnu=0$~~, $\lambda_1\in \mathbb{R}, \lambda_2\in \mathbb{R}.$~~Otherwise,~~$\lambda_1=\lambda_2=0.$~~Case(i) and (ii) are true.
\end{pf}
\vskip 0.5 true cm
\noindent{\bf 3.4 Left-invariant Ricci collineation of $G_4$}\\
\vskip 0.5 true cm
By (2.39) in \cite{C2}, we have
\\
\begin{lem}
Ricci symmetric tensors of $(G_4,g,J)$  associated to Yano connection ~~$\nabla^{*}$~~are given by
\begin{align}
\overline{\rm Ric^{*}}(\overline{e}_i,\overline{e}_j)=\left(\begin{array}{ccc}
-1+n(n-2v)&0&0\\
0&-1+m(n-2v)&\frac{1}{2}m\\
0&\frac{1}{2}m&0
\end{array}\right).
\end{align}
\end{lem}
By (2.2) and (3.16), we have
\begin{lem}
For $G_4$, the following equalities hold
\begin{align}
&{\rm L}_{\overline{e}_1} \overline{\rm Ric^{*}}(\overline{e}_1,\overline{e}_2)=0,~~{\rm L}_{\overline{e}_1}\overline{\rm Ric^{*}}(\overline{e}_1,\overline{e}_3)=0~~,{\rm L}_{\overline{e}_1} \overline{\rm Ric^{*}}(\overline{e}_2,\overline{e}_2)=-2+3m(n-2v),\\\notag
&{\rm L}_{\overline{e}_1} \overline{\rm Ric^{*}}(\overline{e}_2,\overline{e}_3)=n[-1+m(n-2v)]~~,{\rm L}_{\overline{e}_1} \overline{\rm Ric^{*}}(\overline{e}_3,\overline{e}_3)=mn,~~{\rm L}_{\overline{e}_2} \overline{\rm Ric^{*}}(\overline{e}_1,\overline{e}_1)=0~~\\\notag
&{\rm L}_{\overline{e}_2} \overline{\rm Ric^{*}}(\overline{e}_1,\overline{e}_2)=1-\frac{3}{2}m(n-2v),~~{\rm L}_{\overline{e}_2} \overline{\rm Ric^{*}}(\overline{e}_1,\overline{e}_3)=-m[-\frac{1}{2}+n(n-2v)]~~,{\rm L}_{\overline{e}_2} \overline{\rm Ric^{*}}(\overline{e}_2,\overline{e}_3)=0,\\\notag
&{\rm L}_{\overline{e}_2}\overline{\rm Ric^{*}}(\overline{e}_3,\overline{e}_3)=0~~,{\rm L}_{\overline{e}_3} \overline{\rm Ric^{*}}(\overline{e}_1,\overline{e}_1)=0,~~{\rm L}_{\overline{e}_3} \overline{\rm Ric^{*}}(\overline{e}_1,\overline{e}_2)=n-\frac{1}{2}m~~\\\notag
&{\rm L}_{\overline{e}_3}\overline{\rm Ric^{*}}(\overline{e}_1,\overline{e}_3)=-\frac{1}{2}mn,~~{\rm L}_{\overline{e}_3}\overline{\rm Ric^{*}}(\overline{e}_2,\overline{e}_2)=0~~,{\rm L}_{\overline{e}_3} \overline{\rm Ric^{*}}(\overline{e}_2,\overline{e}_3)=0.~~\\\notag
\end{align}
\end{lem}
By (3.16) and definition 2.9, we have
\vskip 0.5 true cm
\begin{thm}
the unimodular Lorentzian Lie group $(G_4,g,J)$ admits left-invariant Ricci collineations associated to the Yano connection $\nabla^{*}$ if and only if\\
\\
$(1)m=0,\lambda_2+n\lambda_3=0,\mathbb{V}_{\mathbb{R}\mathbb{C}}=\langle -n\overline{e}_2+\overline{e}_3\rangle,$\\
\\
$(2)m\neq0,n=0,mv+1=0,\lambda_2=\lambda_3=0,\mathbb{V}_{\mathbb{R}\mathbb{C}}=\langle \overline{e}_1\rangle,$\\
\\
$(3)mn\neq0,\lambda_1=0,v=1,\left|A\right|=0,(2n^2-4n-1)\lambda_2+n\lambda_3=0,\mathbb{V}_{\mathbb{R}\mathbb{C}}=\langle \overline{e}_2-\frac{2n^2-4n-1}{n} \overline{e}_3\rangle,$\\
\\
$(4)mn\neq0,\lambda_1=0,v=-1,\left|A\right|=0,(2n^2+4n-1)\lambda_2+n\lambda_3=0,\mathbb{V}_{\mathbb{R}\mathbb{C}}=\langle \overline{e}_2-\frac{2n^2+4n-1}{n} \overline{e}_3\rangle,$\\
\\
Notice:the value of ~~$\left|A\right|$~~ is in the proof part.
\end{thm}
\begin{pf}
If the unimodular Lorentzian Lie group $(G_4,g,J)$ admits left-invariant Ricci collineations associated to the Yano connection $\nabla^{*}$, then
${\rm L}_{\xi} \overline{\rm Ric^{*}_{ij}}^0=0,i=1,2,3~~ and~~ j=1,2,3$, so
\begin{align}
\left\{\begin{array}{l}
[1-\frac{1}{2}m(n-2v)]\lambda_2+(n-\frac{1}{2}m)\lambda_3=0,\\
\\
m[-\frac{1}{2}+n(n-2v)]\lambda_2+\frac{1}{2}mn\lambda_3=0,\\
\\
\left[-2+3m(n-2v)\right]\lambda_1=0,\\
\\
n[-1+m(n-2v)]\lambda_1=0,\\
\\
mn\lambda_1=0.\\
\end{array}\right.
\end{align}
  We analyze each one of these factors by separate.
  Firs of all, if ~~$m=0$~~, we have
\begin{align}
\left\{\begin{array}{l}
\lambda_2+n\lambda_3=0,\\
\\
\lambda_1=0.\\
\end{array}\right.
\end{align}
So case (1) is true. \\Secondly, when~~$m\neq0,n=0$~~,then
\begin{align}
\left\{\begin{array}{l}
(1+3mv)\lambda_2-\frac{1}{2}m\lambda_3=0,\\
\\
m\lambda_2=0,\\
\\
(1+3mv)\lambda_1=0.\\
\end{array}\right.
\end{align}
From the first and the second equation of the symtem (3.20), there is obvious that ~~$\lambda_2=\lambda_3=0.$~~And as for the third equation, ~~$\lambda_1\in \mathbb{R}$~~ if ~~$1+3mv=0$~~. But if ~~$1+3mv\neq0,~~$,$\lambda_1=0.$~~case (2) is true.\\
Thirdly, if ~~$mn\neq0,$~~then system (3.18)  would reduce to
\begin{align}
\left\{\begin{array}{l}
[1-\frac{3}{2}m(n-2v)]\lambda_2+(n-\frac{1}{2}m)\lambda_3=0,\\
\\
m[-\frac{1}{2}+n(n-2v)]\lambda_2+\frac{1}{2}mn\lambda_3=0,\\
\\
\lambda_1=0.\\
\end{array}\right.
\end{align}
Note that we still have a condition that is not being utillized, which is ~~$v=\pm 1$~~.\\ When ~~$v=1$~~, we assume that
\begin{align}
A=(a_{ij})
,~~\left\{\begin{array}{l}
a_{11}=1-\frac{3}{2}m(n-2),\\
\\
a_{12}=n-\frac{1}{2}m,\\
\\
a_{21}=m[-\frac{1}{2}+n(n-2)],\\
\\
a_{22}=\frac{1}{2}mn.\\
\end{array}\right.~~
\end{align}
We get$$
\left|A\right|
=\frac{1}{4}m[m(n^2-2n+1)+4n(n^2-2n-1)]
$$
If~~$\left|A\right|\neq0$~~, we have ~~$\lambda_2=\lambda_3=0$~~.And if ~~$\left|A\right|=0$~~, we get ~~$(2n^2-4n-1)\lambda_2+n\lambda_3=0.$~~ case (3) is true.\\
When ~~$v=-1$~~, we assume that
\begin{align}
A=(a_{ij})
,~~\left\{\begin{array}{l}
a_{11}=1-\frac{3}{2}m(n+2),\\
\\
a_{12}=n-\frac{1}{2}m,\\
\\
a_{21}=m[-\frac{1}{2}+n(n+2)],\\
\\
a_{22}=\frac{1}{2}mn.\\
\end{array}\right.~~
\end{align}
We get$$
\left|A\right|
=-\frac{1}{4}m[m(n+1)^2+4n(n^2+2n-1)]
$$
If~~$\left|A\right|\neq0$~~, we have ~~$\lambda_2=\lambda_3=0$~~.And if ~~$\left|A\right|\neq0$~~, we have ~~$(2n^2+4n+1)\lambda_2+n\lambda_3=0.$~~Case (4) is true. We complete the proof of theorem.\\
\end{pf}
\vskip 1 true cm
\section{ Left-invariant Ricci collineations associated to canonical connections and Kobayashi-Nomizu connections on three-dimensional non-unimodular Lorentzian Lie groups}

\vskip 0.5 true cm
\noindent{\bf 4.1 Left-invariant Ricci collineation of $G_5$}\\
\vskip 0.5 true cm
By (47) in \cite{C2}, we have
\\
\begin{lem}
Ricci symmetric tensors of $(G_5,g,J)$  associated to Yano connection ~~$\nabla^{*}$~~ are given by
\begin{align}
\overline{\rm Ric^{*}}(\overline{e}_i,\overline{e}_j)=\left(\begin{array}{ccc}
0&0&0\\
0&0&0\\
0&0&0
\end{array}\right).
\end{align}
\end{lem}
By (2.2) and (4.1), we have
\begin{lem}
For $G_5$, the following equalities hold
\begin{align}
&{\rm L}_{\overline{e}_1} \overline{\rm Ric^{*}}(\overline{e}_1,\overline{e}_2)=0,~~{\rm L}_{\overline{e}_1}\overline{\rm Ric^{*}}(\overline{e}_1,\overline{e}_3)=0~~,{\rm L}_{\overline{e}_1} \overline{\rm Ric^{*}}(\overline{e}_2,\overline{e}_2)=0,\\\notag
&{\rm L}_{\overline{e}_1} \overline{\rm Ric^{*}}(\overline{e}_2,\overline{e}_3)=0~~,{\rm L}_{\overline{e}_1} \overline{\rm Ric^{*}}(\overline{e}_3,\overline{e}_3)=0,~~{\rm L}_{\overline{e}_2} \overline{\rm Ric^{*}}(\overline{e}_1,\overline{e}_1)=0~~\\\notag
&{\rm L}_{\overline{e}_2} \overline{\rm Ric^{*}}(\overline{e}_1,\overline{e}_2)=0,~~{\rm L}_{\overline{e}_2} \overline{\rm Ric^{*}}(\overline{e}_1,\overline{e}_3)=0~~,{\rm L}_{\overline{e}_2} \overline{\rm Ric^{*}}(\overline{e}_2,\overline{e}_3)=0,\\\notag
&{\rm L}_{\overline{e}_2}\overline{\rm Ric^{*}}(\overline{e}_3,\overline{e}_3)=0~~,{\rm L}_{\overline{e}_3} \overline{\rm Ric^{*}}(\overline{e}_1,\overline{e}_1)=0,~~{\rm L}_{\overline{e}_3} \overline{\rm Ric^{*}}(\overline{e}_1,\overline{e}_2)=0~~\\\notag
&{\rm L}_{\overline{e}_3}\overline{\rm Ric^{*}}(\overline{e}_1,\overline{e}_3)=0,~~{\rm L}_{\overline{e}_3}\overline{\rm Ric^{*}}(\overline{e}_2,\overline{e}_2)=0~~,{\rm L}_{\overline{e}_3} \overline{\rm Ric^{*}}(\overline{e}_2,\overline{e}_3)=0.~~\\\notag
\end{align}
\end{lem}
By (4.2) and definition 2.9, we have
\vskip 0.5 true cm
\begin{thm}
the non-unimodular Lorentzian Lie group $(G_5,g,J)$ admits left-invariant Ricci collineations associated to the Yano connection $\nabla^{*}$  for ~~$\lambda_1\in \mathbb{R},\lambda_3\in \mathbb{R},\lambda_2\in \mathbb{R}.$~~\\
\end{thm}
\vskip 0.8 true cm
\noindent{\bf 4.2 Left-invariant Ricci collineation of $G_6$}\\
\vskip 0.5 true cm
By (54) in \cite{C2}, we have
\\
\begin{lem}
Ricci symmetric tensors of $(G_6,g,J)$  associated to Yano connection ~~$\nabla^{*}$~~are given by
\begin{align}
\overline{\rm Ric^{*}}(\overline{e}_i,\overline{e}_j)=\left(\begin{array}{ccc}
-(m^2+nu)&0&0\\
0&-m^2&0\\
0&0&0
\end{array}\right).
\end{align}
\end{lem}
By (2.2) and (4.3), we have
\begin{lem}
For $G_6$, the following equalities hold
\begin{align}
&{\rm L}_{\overline{e}_1} \overline{\rm Ric^{*}}(\overline{e}_1,\overline{e}_2)=0,~~{\rm L}_{\overline{e}_1}\overline{\rm Ric^{*}}(\overline{e}_1,\overline{e}_3)=0~~,{\rm L}_{\overline{e}_1} \overline{\rm Ric^{*}}(\overline{e}_2,\overline{e}_2)=2m^3,\\\notag
&{\rm L}_{\overline{e}_1} \overline{\rm Ric^{*}}(\overline{e}_2,\overline{e}_3)=m^2u~~,{\rm L}_{\overline{e}_1} \overline{\rm Ric^{*}}(\overline{e}_3,\overline{e}_3)=0,~~{\rm L}_{\overline{e}_2} \overline{\rm Ric^{*}}(\overline{e}_1,\overline{e}_1)=0~~\\\notag
&{\rm L}_{\overline{e}_2} \overline{\rm Ric^{*}}(\overline{e}_1,\overline{e}_2)=-m^3,~~{\rm L}_{\overline{e}_2} \overline{\rm Ric^{*}}(\overline{e}_1,\overline{e}_3)=0~~,{\rm L}_{\overline{e}_2} \overline{\rm Ric^{*}}(\overline{e}_2,\overline{e}_3)=0,\\\notag
&{\rm L}_{\overline{e}_2}\overline{\rm Ric^{*}}(\overline{e}_3,\overline{e}_3)=0~~,{\rm L}_{\overline{e}_3} \overline{\rm Ric^{*}}(\overline{e}_1,\overline{e}_1)=0,~~{\rm L}_{\overline{e}_3} \overline{\rm Ric^{*}}(\overline{e}_1,\overline{e}_2)=-m^2u~~\\\notag
&{\rm L}_{\overline{e}_3}\overline{\rm Ric^{*}}(\overline{e}_1,\overline{e}_3)=0,~~{\rm L}_{\overline{e}_3}\overline{\rm Ric^{*}}(\overline{e}_2,\overline{e}_2)=0~~,{\rm L}_{\overline{e}_3} \overline{\rm Ric^{*}}(\overline{e}_2,\overline{e}_3)=0.~~\\\notag
\end{align}
\end{lem}
By (4.4) and definition 2.9, we have
\vskip 0.5 true cm
\begin{thm}
the non-unimodular Lorentzian Lie group $(G_6,g,J)$ admits left-invariant Ricci collineations associated to the Yano connection $\nabla^{*}$ if and only if\\
\\
$(1)m=0,v\neq0,\lambda_1\in\mathbb{R},\lambda_2\in\mathbb{R},\lambda_3\in\mathbb{R},\mathbb{V}_{\mathbb{R}\mathbb{C}}=\langle \overline{e}_1,\overline{e}_2,\overline{e}_3\rangle,$\\
\\
$(2)m\neq0,\lambda_1=0,m\lambda_2+u\lambda_3=0,\mathbb{V}_{\mathbb{R}\mathbb{C}}=\langle -\frac{u}{m} \overline{e}_2+ \overline{e}_3\rangle.$\\
\end{thm}
\begin{pf}
If the non-unimodular Lorentzian Lie group $(G_6,g,J)$ admits left-invariant Ricci collineations associated to the Yano connection $\nabla^{*}$, then
${\rm L}_{\xi} \overline{\rm Ric^{*}_{ij}}^0=0,i=1,2,3~~ and~~ j=1,2,3$, so
\begin{align}
\left\{\begin{array}{l}
m^2(m\lambda_2+u\lambda_3)=0,\\
\\
m^3\lambda_1=0,\\
\\
m^2u\lambda_1=0,\\
\\
m+u\neq0,\\
\\
mu-nv=0.\\
\end{array}\right.
\end{align}
  We analyze each one of these factors by separate.
  Firs of all, if ~~$m=0$~~, naturally we have~~$v\neq0,n=0,\lambda_1\in\mathbb{R},\lambda_2\in\mathbb{R},\lambda_3\in\mathbb{R}$ .Case (1) is true.\\
  Then when~~$m\neq0$~~, naturally $\lambda_1$~~must be equal to zero from the second equation and also get ~~$m\lambda_2+u\lambda_3=0$~~.Case (2) is true.
\end{pf}
\vskip 0.5 true cm
\noindent{\bf 4.3 Left-invariant Ricci collineation of $G_7$}\\
\vskip 0.5 true cm
By (64) in \cite{C2} , we have
\\
\begin{lem}
Ricci symmetric tensors of $(G_7,g,J)$  associated to Yano connection ~~$\nabla^{*}$~~ are given by
\begin{align}
\overline{\rm Ric^{*}}(\overline{e}_i,\overline{e}_j)=\left(\begin{array}{ccc}
-m^2&\frac{1}{2}(nv-mn)&mn+nv\\
\frac{1}{2}(nv-mn)&-(m^2+n^2+nu)&\frac{1}{2}(2v^2+nu+mv)\\
mn+nv&\frac{1}{2}(2v^2+nu+mv)&0
\end{array}\right).
\end{align}
\end{lem}
By (2.2) and (4.6), we have
\begin{lem}
For $G_7$, the following equalities hold
\begin{align}
&{\rm L}_{\overline{e}_1} \overline{\rm Ric^{*}}(\overline{e}_1,\overline{e}_2)=-m^3+n^2(\frac{1}{2}m+\frac{3}{2}v),~~{\rm L}_{\overline{e}_1}\overline{\rm Ric^{*}}(\overline{e}_1,\overline{e}_3)=m^3-n^2(\frac{1}{2}m+\frac{3}{2}v)~~,{\rm L}_{\overline{e}_1} \overline{\rm Ric^{*}}(\overline{e}_2,\overline{e}_2)=-3m^2n+2mnv-n^2u-2n^3+2nv^2,\\\notag
&{\rm L}_{\overline{e}_1} \overline{\rm Ric^{*}}(\overline{e}_2,\overline{e}_3)=n(\frac{5}{2}m^2+n^2+nu+\frac{1}{2}mv)~~,{\rm L}_{\overline{e}_1} \overline{\rm Ric^{*}}(\overline{e}_3,\overline{e}_3)=-n(2m^2+3mv+nu+2v^2),~~\\\notag
&{\rm L}_{\overline{e}_2} \overline{\rm Ric^{*}}(\overline{e}_1,\overline{e}_1)=2m^2-n(mn+3nv),~~{\rm L}_{\overline{e}_2} \overline{\rm Ric^{*}}(\overline{e}_1,\overline{e}_2)=n(\frac{3}{2}m^2+n^2+\frac{1}{2}nu-mv-v^2),~~\\\notag
&{\rm L}_{\overline{e}_2} \overline{\rm Ric^{*}}(\overline{e}_1,\overline{e}_3)=-n(\frac{5}{2} v^2+\frac{1}{2}nu+m^2+2mv)+m^2u~~,{\rm L}_{\overline{e}_2} \overline{\rm Ric^{*}}(\overline{e}_2,\overline{e}_3)=\frac{1}{2}mnu+\frac{1}{2}mv^2+n^2v-v^3,\\\notag
&{\rm L}_{\overline{e}_2}\overline{\rm Ric^{*}}(\overline{e}_3,\overline{e}_3)=-2mnu-3nuv-mv^2-2v^3~~,{\rm L}_{\overline{e}_3} \overline{\rm Ric^{*}}(\overline{e}_1,\overline{e}_1)=-2m^3+mn^2+3n^2v,~~\\\notag
&{\rm L}_{\overline{e}_3} \overline{\rm Ric^{*}}(\overline{e}_1,\overline{e}_2)=n(\frac{5}{2}v^2-\frac{3}{2}m^2+\frac{3}{2}mv-n^2-\frac{1}{2}nu)-m^2u~~,{\rm L}_{\overline{e}_3}\overline{\rm Ric^{*}}(\overline{e}_1,\overline{e}_3)=n(m^2+v^2+\frac{3}{2}mv+\frac{1}{2}nu),~~\\\notag
&{\rm L}_{\overline{e}_3}\overline{\rm Ric^{*}}(\overline{e}_2,\overline{e}_2)=-mnu-m^2v-2n^2v+2v^3~~,{\rm L}_{\overline{e}_3} \overline{\rm Ric^{*}}(\overline{e}_2,\overline{e}_3)=mnu+\frac{3}{2}nuv+\frac{1}{2}mv^2+v^3.~~\\\notag
\end{align}
\end{lem}
By (4.7) and definition 2.9, we have
\vskip 0.5 true cm
\begin{thm}
the non-unimodular Lorentzian Lie group $(G_7,g,J)$ admits left-invariant Ricci collineations associated to the Yano connection $\nabla^{*}$ if and only if\\
$(1)mn=0,v\neq0,\lambda_1\in\mathbb{R},\lambda_2=\lambda_3=0,\mathbb{V}_{\mathbb{R}\mathbb{C}}=\langle \overline{e}_1\rangle,$\\
\\
$(2)m=0,v\neq0,u=0,n\lambda_1+v\lambda_2=0,\mathbb{V}_{\mathbb{R}\mathbb{C}}=\langle-\frac{v}{n} \overline{e}_1+\overline{e}_2+\overline{e}_3\rangle$\\
\\
$(3)m\neq0,u=0,n=\sqrt{2}m,v=0,\lambda_1=0,\lambda_2+\lambda_3=0,\mathbb{V}_{\mathbb{R}\mathbb{C}}=\langle \overline{e}_2+\overline{e}_3\rangle,$\\
\\
$(4)m\neq0,u=0,2m^3-n^2(m+3v)\neq0,\lambda_1=0,\lambda_2+\lambda_3=0,\mathbb{V}_{\mathbb{R}\mathbb{C}}=\langle \overline{e}_2+\overline{e}_3\rangle,$\\
\\
$(5)m\neq0,v=0,n=0,\lambda_1=0,\lambda_2=\lambda_3\in \mathbb{R},\mathbb{V}_{\mathbb{R}\mathbb{C}}=\langle \overline{e}_2+\overline{e}_3\rangle,$\\
\\
$(6)m^2n+2mn\neq0,v=0,\lambda_1=0,\lambda_2=\lambda_3\in \mathbb{R},\mathbb{V}_{\mathbb{R}\mathbb{C}}=\langle \overline{e}_2+\overline{e}_3\rangle$\\
\\
$(7)m=-2,v=0,\lambda_1=0,n\neq0,\lambda_2=\lambda_3\in \mathbb{R},\mathbb{V}_{\mathbb{R}\mathbb{C}}=\langle \overline{e}_2+\overline{e}_3\rangle.$\\
\end{thm}
\vskip 0.5 true cm
\begin{pf}
If the non-unimodular Lorentzian Lie group $(G_7,g,J)$ admits left-invariant Ricci collineations associated to the Yano connection $\nabla^{*}$, then
${\rm L}_{\xi} \overline{\rm Ric^{*}_{ij}}^0=0,i=1,2,3~~ and~~ j=1,2,3$, so
\begin{align}
\left\{\begin{array}{l}
[2m^3-n(mn+3nv)](\lambda_2-\lambda_3)=0,\\
\\
\left[-m^3+n^2(\frac{1}{2}m+\frac{3}{2}v)\right]\lambda_1+[n(\frac{3}{2}m^2+n^2+\frac{1}{2}nu-mv-v^2)\lambda_2+[n(\frac{5}{2}v^2-\frac{3}{2}m^2+\frac{3}{2}mv-n^2-\frac{1}{2}nu)-m^2u]\lambda_3=0,\\
\\
\left[m^3-n^2(\frac{1}{2}m+\frac{3}{2}v)\right]\lambda_1+[-n(\frac{5}{2}v^2+m^2+2mv+\frac{1}{2}nu)+m^2u]\lambda_2+n(m^2+v^2+\frac{3}{2}mv+\frac{1}{2}nu)\lambda_3=0,\\
\\
(-3m^2n+2mnv-n^2u-2n^3+2nv^2)\lambda_1+(-m^2v-2n^2v+2v^3)\lambda_3=0,\\
\\
n(\frac{5}{2}m^2+n^2+nu+\frac{1}{2}mv)\lambda_1+(\frac{1}{2}mv^2+n^2v-v^3)\lambda_2+(\frac{3}{2}nuv+\frac{1}{2}mv^2+v^3)\lambda_3=0,\\
\\
n(2m^2+3mv+nu+2v^2)\lambda_1+(3nuv+mv^2+2v^3)\lambda_2=0,\\
\\
mu=0,m+v\neq0.\\
\end{array}\right.
\end{align}
We analyze each one of these factors by separate.
First of all, if ~~$m=0$~~, naturally we have~~$v\neq0$~~and this system (4.8) would reduce to the following system:
\begin{align}
\left\{\begin{array}{l}
n^2v(\lambda_2-\lambda_3)=0,\\
\\
\frac{3}{2}n^2v\lambda_1+n(n^2+\frac{1}{2}nu-v^2)\lambda_2+n(\frac{5}{2}v^2-n^2-\frac{1}{2}nu)\lambda_3=0,\\
\\
\left[\frac{3}{2}n^2v\right]\lambda_1+n(\frac{1}{2}nu+\frac{5}{2}v^2)\lambda_2-n(v^2+\frac{1}{2}nu)\lambda_3=0,\\
\\
(-n^2u-2n^3+2nv^2)\lambda_1+(-2n^2v+2v^3)\lambda_3=0,\\
\\
n(n^2+nu)\lambda_1+(n^2v-v^3)\lambda_2+(\frac{3}{2}nuv+v^3)\lambda_3=0,\\
\\
n(nu+2v^2)\lambda_1+(3nuv+2v^3)\lambda_2=0.\\
\end{array}\right.
\end{align}
Notice that if~~$n=0$~~at same time, the system (4.9) reduces to
\begin{align}
\left\{\begin{array}{l}
v^3\lambda_3=0,\\
\\
v^3(\lambda_2-\lambda_3)=0,\\
\\
v^3\lambda_2=0.\\
\end{array}\right.
\end{align}
Because~~$v\neq0$~~, so ~~$\lambda_2=\lambda_3=0$~~and ~~$\lambda_1\in\mathbb{R}$~~.Case (1) is true.\\
As for~~$n\neq0$~~, we would get~~$\lambda_2=\lambda_3$~~from the first equation. Put~~$\lambda_2=\lambda_3$~~into the second and sixth equation of the system (4.9), we have~~$nuv\lambda_2=0$~~. Then if~~$u=0$~~, we have~~$\lambda_2\in \mathbb{R},\lambda_3\in \mathbb{R}, n\lambda_1+v\lambda_2=0$~~. put ~~$n\lambda_1+v\lambda_2=0$~~ into other equations of the system (4.9), they did not produce contradictions. So case (2) is true.\\
Next when~~$m\neq0$~~ and ~~$ u=0$,~~ system (4.8) would reduce to
\begin{align}
\left\{\begin{array}{l}
[2m^3-n^2(m+3v)](\lambda_2-\lambda_3)=0,\\
\\
-\frac{1}{2}[2m^3-n^2(m+3v)]\lambda_1+n(\frac{3}{2}m^2+n^2-mv-v^2)\lambda_2+n(\frac{5}{2}v^2-\frac{3}{2}m^2+\frac{3}{2}mv-n^2)\lambda_3=0,\\
\\
\frac{1}{2}[2m^3-n^2(m+3v)]\lambda_1-n(\frac{5}{2}v^2+m^2+2mv)\lambda_2+n(m^2+v^2+\frac{3}{2}mv)\lambda_3=0,\\
\\
(-3m^2n+2mnv-2n^3+2nv^2)\lambda_1+(-m^2v-n^2v+2v^3)\lambda_3=0,\\
\\
n(\frac{5}{2}m^2+n^2+\frac{1}{2}mv)\lambda_1+\frac{1}{2}mv^2+n^2v-v^3)\lambda_2+(\frac{1}{2}mv^2+v^3)\lambda_3=0,\\
\\
n(2m^2+3mv+2v^2)\lambda_1+(mv^2+2v^3)\lambda_2=0.\\
\end{array}\right.
\end{align}
If ~~$\left[2m^3-n(mn+3nv)\right]=0$~~, notice that~~$n\neq 0$~~ or contradiction. First of all, let's consider the first two equations temporarily:
\begin{align}
\left\{\begin{array}{l}
(\frac{3}{2}m^2+n^2-mv-v^2)\lambda_2+(\frac{5}{2}v^2-\frac{3}{2}m^2+\frac{3}{2}mv-n^2)\lambda_3=0,\\
\\
(\frac{5}{2}v^2+m^2+2mv)\lambda_2-(m^2+v^2+\frac{3}{2}mv)\lambda_3=0.\\
\\
2m^3-n^2(m+3v)=0.\\
\end{array}\right.
\end{align}
we assume that
\begin{align}
A=(a_{ij}),~~\left\{\begin{array}{l}
a_{11}=\frac{3}{2}m^2+n^2-mv-v^2,\\
\\
a_{12}=\frac{5}{2}v^2-\frac{3}{2}m^2+\frac{3}{2}mv-n^2,\\
\\
a_{21}=\frac{5}{2}v^2+m^2_2mv,\\
\\
a_{22}=-(m^2+v^2+\frac{3}{2}mv).\\
\end{array}\right.
\end{align}
We get~~$\left|A\right|
=-\frac{1}{4}v((-m^3+3m^2v-2mn^2+25mv^2-6n^2v+21v^3)
$~~and combine this with the above equation of the system (4.13), have the following two solutions:
\begin{align}
\left\{\begin{array}{l}
\left|A\right|=0,n=-\sqrt{2}m,v=0,\lambda_1=0,\lambda_2=\lambda_3,\\
\\
\left|A\right|=0,n=\sqrt{2}m,v=0,\lambda_1=0,\lambda_2=\lambda_3.\\
\end{array}\right.
\end{align}
Case (3)and Case (4) are true.But if~~ $\left[2m^3-n(mn+3nv)\right]\neq0$~~,naturally we have~~$\lambda_2=\lambda_3$.  Put this into system (4.11) and get
\begin{align}
\left\{\begin{array}{l}
\frac{1}{2}\left[2m^3-n(mn+3nv)\right]+\lambda_1+n(\frac{1}{2}v^2+\frac{3}{2}mv)\lambda_2,\\
\\
n(\frac{5}{2}m^2+n^2+\frac{1}{2}mv)\lambda_1+(m^2v+n^2v)\lambda_2=0,\\
\\
(-3m^2n+2mnv-2n^3+2nv^2)\lambda_1+(2v^3-m^2v-n^2v)\lambda_2=0,\\
\\
n(2m^2+3mv+2v^2)\lambda_1+(mv^2+2v^3)\lambda_2=0.\\
\end{array}\right.
\end{align}
Since the number of unknowns is less than the number of equations, we first consider the second and third equations.we assume that
\begin{align}
B=(b_{ij}),~~\left\{\begin{array}{l}
b_{11}=n(\frac{5}{2}m^2+n^2+\frac{1}{2}mv),\\
\\
b_{12}=mv^2+n^2v,\\
\\
b_{21}=n(2m^2+3mv+2v^2),\\
\\
b_{22}=mv^2+2v^3.\\
\end{array}\right.
\end{align}
We get~~$\left|B\right|
=\frac{1}{2}nv((m^3v-4m^2n^2-5m^2v^2+10m^2v-4mn^2v-4mv^3+2mv^2-4n^2v^2+4n^2v)
$~~and if $\left|B\right|\neq0$~~combining this with the above equation of the system (4.15), have the following solutions:
\begin{align}
\left\{\begin{array}{l}
m\neq0,v=0,n=0,\lambda_1=0,\lambda_2=\lambda_3\in \mathbb{R},\\
\\
m^2n+2mn\neq0,v=0,\lambda_1=0,\lambda_2=\lambda_3\in \mathbb{R},\\
\\
m=-2,v=0,\lambda_1=0,n\neq0,\lambda_2=\lambda_3\in \mathbb{R},\\
\\
m\neq0,v\neq0,n=0,\lambda_2=\lambda_3=0,\lambda_1\in\mathbb{R}\\
\end{array}\right.
\end{align}
We substitute these four results into the remaining equations for testing. After testing, we found that the first three results do not conflict. Case (5), Case (6) and Case(7) are ture, thus we complete the proof of the theorem.
\end{pf}
\vskip 1 true cm
\section{Acknowledgements}

The author are deeply grateful to the referees for their valuable commments and helpful suggestions.
\vskip 1 true cm


\bigskip
\bigskip

\noindent {\footnotesize {\it Yu Tao} \\
{School of Mathematics and Statistics, Northeast Normal University, Changchun 130024, China}\\
{Email: yut338@nenu.edu.cn}\\

\end{document}